\documentclass[english,a4paper,12pt]{article}
\usepackage{babel}
\usepackage[ansinew]{inputenc}
\usepackage{graphicx}
\usepackage{verbatim}
\usepackage{amsmath}
\usepackage{amssymb}
\def\C{\mathbb C}

\def\Z{\mathbb Z}

\newtheorem{theo}{Theorem}[section]
\newtheorem{ex}[theo]{Example}
\newtheorem{lemma}[theo]{Lemma}
\newtheorem{prop}[theo]{Proposition}
\newtheorem{defi}[theo]{Definition}
\newtheorem{rem}[theo]{Remark}

\def\pr{{\it Proof. }}

\begin{document}

\title{A uniqueness problem concerning entire functions and their derivatives}
\author{Andreas Sauer and Andreas Schweizer}

\date{\today}

\maketitle

\begin{abstract}
We determine all entire functions $f$, such that for nonzero complex values $a \neq b$ the implications $f=a \Rightarrow f'=a$ and $f'=b \Rightarrow f=b$ hold. This solves an open problem in uniqueness theory. In this context we give a normality criterion, which might be interesting in its own right.
\end{abstract}

\noindent
\section{Introduction}
Two meromorphic functions $f$ and $g$ on a complex domain $D$ are 
said to share the value $a\in\C$ if for all $z_0\in D$ we have
$f(z_0)=a\Leftrightarrow g(z_0)=a$.
Usually this is more loosely written as 
$f=a\Leftrightarrow g=a$.
Often one also says more precisely that $f$ and $g$ are sharing 
the value $a$ IM (ignoring multiplicities) to emphasize that one
does not require that at such a point the function $f$ takes the
value $a$ with the same multiplicity as $g$.

Here we are mainly interested in the case where the domain $D$ is 
the whole complex plane $\C$, $f$ is an entire function, and $g$
is its derivative $f'$. The most famous theorem in this setting 
is the following result from 1979.
\\ \\
{\bf Theorem A.} \cite[Satz 1]{MuSt} \it
Let $f$ be a nonconstant entire function and let $a$ and $b$ be complex 
numbers with $a\neq b$. If $f$ and $f'$ share the values $a$ and $b$ IM, 
then $f\equiv f'$.
\rm
\\ \\
See also \cite[Theorem 8.3]{YaYi} for the same proof in English. A different 
proof for the case $ab\neq 0$ was also given in \cite[Theorem 2]{Gu}.
\par
Theorem A has been generalized in many ways. Specifically we mention the f
ollowing result which gives up one of the directions of the implication 
$f=b\Leftrightarrow f'=b$. Actually, it also provides a third different 
proof for the case $ab\neq 0$ of Theorem A.
\\ \\
{\bf Theorem B.} \cite[Theorem 1.1]{LuXuYi} \it
Let $a$ and $b$ be two nonzero distinct complex numbers, 
and let $f(z)$ be a nonconstant entire function. If 
$f(z)=a\Leftrightarrow f'(z)=a$ and 
$f'(z)=b\Rightarrow f(z)=b$, then $f\equiv f'$.
\rm
\\ \\
The functions $f=\frac{b}{4}z^2$ and $f=C\text{\rm e}^{\frac{a}{C}z}+a-C$
from \cite{LiYi} show that the conditions $a\neq 0$ and $b\neq 0$
cannot be removed in Theorem B. Here $C$ is any nonzero constant.

In \cite[Question 2]{LuXuYi} it is asked what happens if in Theorem B the condition 
$f(z)=a \Leftrightarrow f'(z)=a$ is weakened to
$f(z)=a \Rightarrow f'(z)=a$. 
We will answer this question by proving the following theorem, which is the main result of this paper.

\begin{theo} \label{main} {\rm (Main Theorem)} \\
Let $f$ be a nonconstant entire function and let $a\neq b$ be two nonzero complex numbers. If
$$f(z_0)=a \Rightarrow f'(z_0)=a$$
and
$$f'(z_0)=b \Rightarrow f(z_0)=b$$
for all $z_0 \in \C$, then $f$ is one of the following:
\begin{itemize}
\item[(i)] $f(z)=az+C$ where $C$ is any constant. These are all nonconstant polynomial solutions.
\item[(ii)] $f\equiv f'$, i.e. $f(z)=C\text{\rm e}^z$ where $C$ is any nonzero constant.
\item[(iii)] $f(z)=C\text{\rm e}^{\frac{b}{b-a}z} +a$ where $C$ is any nonzero constant. In this case $f$ and $f'$ do even share the value $b$. Moreover, $a$ is a Picard value of $f$.
\item[(iv)] $f(z)=6aC\text{\rm e}^{\frac{z}{6}}(C\text{\rm e}^{\frac{z}{6}}-1)+a$ with $b=\frac{-a}{8}$ and $C$ is any nonzero constant.
\end{itemize}
\end{theo}

The condition $f=a \Rightarrow f'=a$ holds of course trivially if $a$ is a Picard value of $f$, i.e.~if $f$ omits the value $a$. This would give a quick corollary to our main theorem. But actually we have to go the other way round and first prove the following more general version of this corollary, because we need it as one of the many steps in the long proof of our main theorem.

\begin{theo} \label{b=>b} Let $f$ be a transcendental entire function with
$f'=b\Rightarrow f=b$ for some nonzero value $b$.
If the value $a$ is taken by $f$ only finitely often, 
then $a\neq b$ and
$$f(z)=C{\rm e}^{\frac{b}{b-a}z}+a$$
where $C$ is a nonzero constant. In particular, $a$ is a Picard exceptional value of $f$, all $b$-points of $f'$ are simple and $b$ is a shared value of $f$ and $f'$. Furthermore,
if $a=0$ we have $f\equiv f'$. 
\end{theo}

In Section \ref{sectionnormal} we will prove Theorem \ref{b=>b} after first showing a normality criterion that we will need for both, Theorem \ref{main} and \ref{b=>b}.

The proof of Theorem \ref{main} is much longer and requires lots of different arguments, some even from
Number Theory. Therefore we now give some sort of rough overview. 
\par
Using the normality criterion from Section \ref{sectionnormal} we can show that $f$ satisfies a differential
equation 
$$\frac{f''(f' -f)}{(f-a)(f'-b)}=c$$
with a constant $c$. If $c=0$ we obtain the solutions $(i)$ and $(ii)$ in Theorem \ref{main}.
All other solutions must be transcendental (for degree reasons). If such a solution $f$ (for $c\neq 0$) 
or one of its derivatives $f^{(n)}$ takes a value from $\C$ only finitely often,  we can show with the help of 
Theorem \ref{b=>b} that $f$ must be of the form $(iii)$ in Theorem \ref{main}.
\par
The next key step is to show that any transcendental solution must be periodic (Theorem \ref{f_periodic}),
and even a polynomial in $\text{\rm e}^{\lambda z}$ where $\frac{2\pi i}{\pm \lambda}$ is a fundamental period with sign chosen
suitably (Lemma \ref{P(exp)}).
With this the counting functions from Nevanlinna Theory for $f$ translate into how often the polynomial
$P$ takes certain values (Definition \ref{d_j_k} and Lemma \ref{countingfct}).
\par
Moreover, using the differential equation, the $n$-th derivative $f^{(n)}$ of $f$ at a point (more precisely, 
a factor $C_n$ times this derivative) can be expressed as a formula in the lower derivatives of $f$ at 
this point. So for example the Taylor expansion of $f$ at an $a$-point of $f$ is essentially determined by
the first few coefficients.This furnishes an upper bound for the degree of the polynomial $P$. See for 
example Lemma \ref{d=4_j=2_k=2}.
The crux is the condition that $C_n$ does not vanish for any $n$. The most natural way to prove this 
non-vanishing (and actually the only one we see) is with arguments from Number Theory. See Lemmas
\ref{k_n_square} to \ref{k=3}.
\par
Once one has a good bound for the degree of $P$, one can take a general form of $P$ and an 
undetermined $\lambda$ and check for which specific values they solve the uniqueness problem.
This is what we do at the end of Sections \ref{sectionbpoints} and \ref{sectionendofproof}.

\section{A normality criterion}\label{sectionnormal}

We will now prove a normality criterion which is a generalization of a criterion given by 
Li and Yi \cite[Theorem 3]{LiYi}. It is slightly more general than we will actually need for 
our purposes, but our Theorem might be interesting in its own right.

\begin{theo} \label{normality} Let ${\cal F}$ be a family of holomorphic functions in a domain 
$D$. Let $a$ and $b \neq 0$ be two complex numbers such that
$b\neq a$, and let $S_a$ and $S_b$ be compact subsets of $\C$ 
with $b\not\in S_a$ and $a\not\in S_b$. If, for each $f\in {\cal F}$ and $z\in D$,
$f=a\Rightarrow f'\in S_a$, and $f'=b\Rightarrow f \in S_b$, 
then ${\cal F}$ is normal in $D$.
\end{theo}

\noindent
\pr The special case $S_a=\{a\}$ and $S_b=\{b\}$ is Theorem 3 from \cite{LiYi}, and our 
proof follows their proof closely, while simplifying some of the arguments.

Setting $h(z)=f(z)-a$ we have $|h'(z)|\leq M + 1$ when $h(z)=0$
where $M$ is the maximum of the absolute values of the elements
of $S_a$. We can assume that $D$ is the unit disk and that the
family $\{f(z)-a\ :\ f\in {\cal F} \}$ is not normal at $0$.

To bring this to a contradiction, as in \cite{LiYi}, we use a version
of the Pang-Zalcman Lemma \cite[Lemma 2]{PaZa}, where
$$g_n(\xi)=\rho_n^{-1}\{f_n(z_n+\rho_n\xi)-a\}\to g(\xi)$$
locally uniformly on $\C$, with a nonconstant entire function $g(\xi)$ satisfying
$$g^\#(\xi)\leq g^\#(0)=M+2.$$
It is immediately clear that $f_n(z_n+\rho_n\xi)$ converges locally uniformly in $\C$ to the constant $a$, since otherwise because of $\rho_n^{-1} \to \infty$ there exists $\xi$ with $g_n(\xi) \to \infty$, which is impossible.

Now we prove $g=0\Rightarrow g'\in S_a$. If $g(\xi_0)=0$, by
Hurwitz's Theorem there exist $\xi_n\to\xi_0$ as $n\to\infty$
such that for sufficiently large $n$ 
$$g_n(\xi_n)=\rho_n^{-1}\{f_n(z_n+\rho_n\xi_n)-a\}=0.$$
So $f_n(z_n+\rho_n\xi_n)=a$ and hence
$g_n'(\xi_n)=f_n'(z_n+\rho_n\xi_n)\in S_a$. Since 
$g_n'(\xi_n)\to g'(\xi_0)$ and $S_a$ is compact, we get 
$g'(\xi_0)\in S_a$.

Next we prove that $g'(\xi)\neq b$ on $\C$. If 
$g'(\xi)\equiv b$, then $g(\xi)=b\xi +c$, which together with 
$g=0\Rightarrow g'\in S_a$ contradicts $b\not\in S_a$. If
$g'(\xi)\not\equiv b$ and $g'(\xi_0) = b$, then again by Hurwitz's Theorem there exist $\xi_n\to\xi_0$ as $n\to\infty$
such that for sufficiently large $n$ 
$$g_n'(\xi_n)=f_n'(z_n+\rho_n\xi_n)=b.$$
By assumption we get $f_n(z_n+\rho_n\xi_n) \in S_b$. Because $S_b$ is compact and $a \not \in S_b$ we get the contradiction $f_n(z_n+\rho_n\xi_n) \not \to a$.

Since $g'(\xi)$ is of exponential type, we now have
$g'(\xi)=b+{\text e}^{b_0 +b_1 \xi}$ with constants $b_0,b_1$ and therefore 
$g(\xi)=c + b\xi + {\text e}^{b_0 +b_1 \xi}/b_1$.
\par

If $b_1\neq 0$ then since $b \neq 0$ it follows that $c + b\xi \not \equiv 0$ is a non-deficient 
small function for ${\text e}^{b_0 +b_1 \xi}/b_1$. Hence there is a sequence $\xi_n \to \infty$ 
with $g(\xi_n)=0$. Since $c + b \xi_n \to \infty$ it follows immediately that 
${\text e}^{b_0 +b_1 \xi_n}/b_1 \to \infty$. This implies that 
$g'(\xi_n)=b+{\text e}^{b_0 +b_1 \xi_n} \to \infty$, but this contradicts 
$g=0\Rightarrow g'\in S_a$ since $S_a$ is compact. 
\par
If $b_1 =0$, we have $g(\xi)=(b+{\text e}^{b_0})\xi +c$. From 
$g=0\Rightarrow g'\in S_a$ we get $g' = b+{\text e}^{b_0}\in S_a$, and
hence the contradiction $g^{\#}(0) \leq |b+{\text e}^{b_0}| < M+2$. 
$\square$
\vskip.3cm
As usual we can draw immediate consequences for certain entire functions.

\begin{theo} \label{exponential_type} Let $f$ be an entire function and let $a$ and $b \neq 0$ be two finite complex numbers such that
$b\neq a$, and let $S_a$ and $S_b$ be compact subsets of $\C$ 
with $b\not\in S_a$ and $a\not\in S_b$. If for every $z \in \C$,
$f(z)=a\Rightarrow f'(z)\in S_a$, and $f'(z)=b\Rightarrow f(z) \in S_b$, 
then $f^\#$ is bounded on $\C$ and therefore $f$ is of exponential type.
\end{theo}

\noindent
\pr Simply apply Theorem \ref{normality} to the family ${\cal F} := \{ f(z + w) : w \in \C \}$ and note that $f$ is of exponential type by a well known theorem of Clunie and Hayman \cite[Theorem 2]{ClHa}. $\square$

\begin{ex} \rm
The function $f(z)=\exp(z^n) -1$ has order $n$. Since $f=-1 \Rightarrow f'=-1$ and $f'=0 \Rightarrow f=0$ both hold, this shows that the condition $b\neq 0$ cannot be omitted in Theorems \ref{normality} and \ref{exponential_type}.
\end{ex}

\noindent
{\it Proof of Theorem \ref{b=>b}.} If $a$ and $b$ were equal, $f'=b\Rightarrow f=b$ would imply that 
$f$ and $f'$ both take the nonzero value $b$ only finitely often, in contradiction to a famous theorem 
by Hayman \cite[Theorem 3]{Ha}.
\par
Set $S_a := \{f'(z) : f(z) = a \}$ and $S_b := \{ b\}$. Then Theorem \ref{exponential_type} can be applied, and it follows that $f$ is an entire function of exponential type. Hence $f$ has the form $f(z) = P(z)\exp(\lambda z) + a$
with $\lambda \neq 0$ and a polynomial $P$. Thus $f'(z) = (P'(z) + \lambda P(z))\exp(\lambda z)$. Obviously $f'$ has only finitely many zeros and thus there are infinitely many $z_0$ with $f'(z_0) =b$. For these we must have 
$\exp(\lambda z_0) = b/(P'(z_0) + \lambda P(z_0))$. Moreover, 
$b - a = f(z_0) - a = \frac{b P(z_0)}{P'(z_0)+\lambda P(z_0)}$ shows that the
rational function 
$$R(z) = \frac{b P(z)}{P'(z)+\lambda P(z)}$$
takes the value $b - a$ infinitely often. So necessarily $R(z) \equiv b - a$. This shows that
$P$ must be a constant and $\lambda = b/(b-a)$. $\square$

\section{A uniqueness problem}

We now define the uniqueness problem which we will solve. 

\begin{defi} \rm Let $f$ be an entire function. We say that $f$ solves the problem (P), if there exist nonzero complex numbers $a \neq b$ such that for all $z \in \C$
\begin{align*}
f(z)=a \Rightarrow f'(z)=a \quad \text{and} \quad f'(z)=b \Rightarrow f(z)=b. \tag{P}
\end{align*}
\end{defi}

We start by solving (P) in certain simple situations. We will need these results later or simply for completeness.

\begin{lemma} \label{poly} If $f$ is a polynomial that solves (P), then $f = az + B$ with arbitrary $B$ or $f = B$ with $B \neq a$.
\end{lemma}

\noindent
\pr This is completely elementary. A polynomial $f$ of degree $d>1$ cannot satisfy $f=a\Rightarrow f'=a$, because this condition would imply that the $a$-points of $f$ are simple, and so $f$ would have $d$ different $a$-points, but $f'$ cannot have more than $d-1$. The rest is easy. $\square$

\begin{lemma} \label{f_finitely_many_A} Suppose $f$ is a transcendental entire function that solves (P). If $f$ takes a value $A \in \C$ only finitely often, then $f(z) = C\exp(bz/(b-A))+A$ with $C \neq 0$ and $A=a$ or $A=0$. In particular, $f$ has no $A$-points at all.
\end{lemma}

\noindent
\pr By Theorem \ref{b=>b} we only have to show $A=a$ or $A=0$. Assume $A\neq a$. Then by Picard's Theorem there are (infinitely many) points $z_0$ with $f(z_0)=a$. So $C\exp \left( \frac{b}{b-A}z_0 \right)=a-A$. By (P) we also have 
$a=f'(z_0)= \frac{b}{b-A}C\exp \left( \frac{b}{b-A}z_0 \right)=
\frac{b}{b-A}(a-A)$. So $ab-aA=ba-bA$, which implies $A=0$.
$\square$

\begin{lemma} \label{derivative_finitely_many_B} Suppose $f$ is a transcendental entire function that solves (P). If some derivative $f^{(n)}$ takes a value $B \in \C$ only finitely often, then $B=0$ and $f(z) = C\exp(bz/(b-A))+A$ with $C \neq 0$ and $A=a$ or $A=0$. In particular, all derivatives have no zeros at all.
\end{lemma}

\noindent
\pr Suppose $f^{(n)}$ takes the value $B$ only finitely often. Since $f^{(n)}$ is an entire function of exponential type we get $f^{(n)}(z) = P(z) \exp(\lambda z) + B$ with a polynomial $P$. Integrating $n$ times we obtain $f(z) = U(z) \exp(\lambda z) + V(z)$ with polynomials $U$ and $V$. If $f$ has only finitely many $a$-points, then the assertion follows immediately from Lemma \ref{f_finitely_many_A}. Therefore we can assume that there are infinitely many $z_0$ with $f(z_0)=a$ and $U(z_0) \neq 0$. For these we get $\exp(\lambda z_0) = (a - V(z_0)) / U(z_0)$ and $a = f'(z_0) = (U'(z_0) + \lambda U(z_0))(a - V(z_0)) / U(z_0) + V'(z_0)$, so
$$
aU(z_0) = (U'(z_0) + \lambda U(z_0))(a - V(z_0)) + U(z_0) V'(z_0).
$$
The degrees on both sides of the equation have to be equal. Except for the trivial case $V \equiv 0$ it follows that $U$ and $U \cdot V$ have the same degree, which implies that $V \equiv A$ is constant. It follows that $f$ has only finitely many $A$-points and the assertion again follows from Lemma \ref{f_finitely_many_A}. 
$\square$

\vskip.3cm
\noindent
Thus we have identified the cases (i) to (iii) in Theorem \ref{main}. The rest of the paper is devoted to the fact, that the only other solutions of (P) are the ones described in (iv) of Theorem \ref{main}.
\vskip.3cm
Of basic importance is an auxiliary function that was also used in \cite{LuXuYi}. 

\begin{theo}
If $f$ solves (P) then $f$ also is a solution of a differential equation
\begin{equation} \label{g}
g=\frac{f''(f'-f)}{(f-a)(f'-b)}=c \tag{D}
\end{equation}
with a constant $c$. We use $g$ for the auxiliary function on the left.
\end{theo}

\noindent
\pr From the conditions in (P) one easily sees that $g$ is an entire function. Moreover,
$$
g=\frac{f'}{f-a} \frac{f''}{f'-b} 
- \frac{f''}{f'-b} - \frac{f'}{f-a} \frac{a}{b} 
\left( \frac{f''}{f'-b} - \frac{f''}{f'} \right).
$$
As we know from Theorem \ref{exponential_type} that $f$ is an entire function of order at most one, by a Theorem of Ngoan and Ostrovskii (see e.g.~\cite[Theorem 3.5.1]{ChYe}) this implies
$T(r,g)=m(r,g)=o(\log r)$, which means that $g$ is a constant.
$\square$

\begin{rem}\label{PversusD} \rm 
We hasten to point out that although we will determine all solutions of the uniqueness problem (P), 
this is not guaranteed to give all solutions of the differential equation (D). The reason for this is that 
(P) is (potentially) a stronger condition than (D) in the sense that although (P) implies (D), conversely
(D) with $c\neq 0$ only implies the weaker conditions
$$f'(z_0)=b\Rightarrow f(z_0)=b$$
and 
$$f(z_0)=a\Rightarrow(f'(z_0)=a\ \ \hbox{\rm or}\ \ (f''(z_0)=0\ \ \hbox{\rm and}\ \ f'(z_0)\neq 0)).$$
As in our proofs we frequently use the condition $f=a\Rightarrow f'=a$, we might be implicitly suppressing
some solutions of (D) that are not solutions of (P).
\end{rem}




\vskip.3cm
\noindent
From Lemma \ref{derivative_finitely_many_B} we know the solutions of (P), when $f''$ has no zeros. Therefore we start to examin solutions where $f''$ has zeros.

\begin{lemma} \label{g=(k+1)b/(a-b)} Let $f$ solve (P). If $f'' \not \equiv 0$ has a zero of order $k \ge 1$ at $z_0$, then 
\begin{itemize}
\item[a)] $f(z_0) = f'(z_0) = b$,
\item[b)] $g \equiv (k+1)b/(a-b)$ with $g$ from (\ref{g}),
\item[c)] all zeros of $f''$ are of order $k$.
\end{itemize}
\end{lemma}

\noindent
\pr a) If $g \equiv 0$ with $g$ from (\ref{g}), then $f'\equiv f$ or $f'' \equiv 0$. In both cases $f''$ either has no zeros or $f'' \equiv 0$. Hence $g \not \equiv 0$. Therefore $(f(z_0)-a)(f'(z_0)-b)$ has to be zero. All $a$-points of $f$ are simple. So if $f(z_0)=a$ then $(f(z_0)-a)(f'(z_0)-b)$ has a simple zero, while the numerator of $g$ has at least a double zero. This contradicts $g \not \equiv 0$. Hence $f'(z_0)=b$ and by assumption $f(z_0)=b$.\\
b) Locally $f''(z) = A(z-z_0)^k + O((z-z_0)^{k+1})$ with $A \neq 0$ and thus 
$$
f'(z) = b + \frac{A}{k+1}(z-z_0)^{k+1} + O\left((z-z_0)^{k+2}\right)
$$ and 
$$
f(z) = b + b(z-z_0) + \frac{A}{(k+1)(k+2)}(z-z_0)^{k+2} + O \left((z-z_0)^{k+3}\right).
$$
Plugging this into $g$ from (\ref{g}) and letting $z \to z_0$ proves the assertion.\\
c) follows immediately from b). 
$\square$

\section{Periodicity}

One of the main steps in our argumentation is to show that solutions of (P), apart from the trivial polynomial solutions, are periodic. We start by showing that at an $a$-point of $f$ there are only two possible initial conditions for the differential equation given in (\ref{g}).

\begin{lemma} \label{f''_at_a,a} Let $f$ solve (P). If $f(z_0)=a$ and $f''$ has zeros of order $k \ge 1$ (elsewhere) then
\begin{align*} 
f''(z_0) = \frac{a}{2} \pm \frac{\sqrt{a^2 + 4(k+1)ab}}{2}
\end{align*}
\end{lemma}

\noindent
\pr If $f(z_0)=a$ then by assumption $f'(z_0)=a$. Therefore $(f-a)(f'-b)$ has a simple zero at $z_0$. Since $f(z_0) - f'(z_0) = 0$ and $g \not \equiv 0$ it follows $f''(z_0) \neq 0$. We get from Lemma \ref{g=(k+1)b/(a-b)}, b)
$$
\frac{f'-f}{f-a} = \frac{(k+1)b}{a-b} \frac{f'-b}{f''}.
$$
Taking the limit for $z \to z_0$ on both sides, we obtain
$$
\frac{f''(z_0)- a}{a} = \frac{(k+1)b}{f''(z_0)}.
$$
Solving this quadratic equation for $f''(z_0)$ proves the assertion. 
$\square$

\begin{lemma} \label{f_periodic_taylor} An entire function $f$ is periodic if and only if there exist two distinct points $z_1$ and $z_2$, such that the sequences of Taylor coefficients at $z_1$ and $z_2$ are identical, in which case $z_1 - z_2$ is a period of $f$.
\end{lemma}

\noindent
\pr If $f$ is periodic then the existence of $z_1$ and $z_2$ is trivial. If $z_1$ and $z_2$ exist, then we get $f(z + z_1) = f(z + z_2)$ and by the identity theorem $f(z) = f(z + z_1 - z_2)$. 
$\square$
\vskip.3cm

Now we are in the position to show periodicity for all transcendental solutions of (P).

\begin{theo} \label{f_periodic} If $f$ is a transcendental entire function that solves (P) then $f$ is periodic.
\end{theo}

\noindent
\pr Let $z_1$, $z_2$ and $z_3$ be three distinct points such that $f(z_j) = a$ and $f'(z_j)=a$ for $j=1,2,3$. Since by Lemma \ref{f''_at_a,a} there are only two
different possibilities for $f''(z_j)$ we may assume that $f''(z_1) = f''(z_2)$. The strategy of the proof is to show that the Taylor coefficients of $f$ in $z_1$ and $z_2$ are identical, so that by Lemma \ref{f_periodic_taylor} $f$ has to be periodic. This strategy fails in certain exceptional cases, which will be ruled out separately.\\
The first possible exception is that $f$ does not possess three $a$ points. But then Lemma \ref{f_finitely_many_A} immediately shows that $f$ is periodic.\\
By Lemma \ref{derivative_finitely_many_B} we may assume that $f''$ has infinitely many zeros. By Lemma \ref{g=(k+1)b/(a-b)} these zeros are all of order $k$ and $g \equiv (k+1)b/(a-b)$ with $g$ from (\ref{g}). This implies
$$
f''(f'-f) = \frac{(k+1)b}{a-b} (f-a)(f'-b).
$$
Taking the $n$-th derivative on both sides according to the Leibniz rule, we obtain
\begin{align*}
f^{(n+2)}(f'-f) + {n \choose 1} f^{(n+1)}(f''-f') + \ldots + f''(f^{(n+1)} - f^{(n)}) & \\
= \frac{(k+1)b}{a-b}\left[f^{(n)}(f'-b) + \ldots + (f-a) f^{(n+1)} \right] &.
\end{align*}
At $z=z_j$ the highest derivative $f^{(n+2)}$ disappears, since $f'(z_j)-f(z_j)=0$. The second highest derivative $f^{(n+1)}$ appears in three terms: in the second and the last term on the left, and in the last term on the right. The latter is again zero because $f(z_j) - a = 0$. We get
\begin{align} \label{C_{n+1}}
C_{n+1}:= (n+1)f'' - nf'
\end{align}
as the coefficient of $f^{(n+1)}$. So the $(n+1)$th Taylor coefficient can be determined uniquely, depending on the values of $a$, $b$ and $k$, but not on the values $z_j$, if $C_{n+1} \neq 0$ for $n \ge 2$.\\
The rest of the proof deals with the last possible exceptional case: Suppose $C_{n+1}=0$ for some $n \ge 2$ at each $a$-point with at most two exceptions. According to Lemma \ref{f''_at_a,a} the equation $C_{n+1}=0$ can be written as
$$
(n+1) \left( \frac{a}{2} \pm \frac{\sqrt{a^2 + 4(k+1)ab}}{2} \right) - na = 0.
$$
It follows by direct calculation
\begin{align} \label{relation_a_b}
a = \frac{-(n+1)^2 (k+1)}{n} \cdot b.
\end{align}
Since $a$, $b$ and $k$ are given, and since $(n+1)^2/n$ is increasing for $n \ge 1$, we get that $n$ is uniquely determined and therefore has the same value at each $a$-point.\\
We consider
$$
h(z) := \frac{(n+1)f''(z) - nf'(z)}{f(z)-a}.
$$
$h$ is meromorphic with at most two poles at the possible exceptional $a$-points and $m(r,h) = O(\log(r))$. Hence $h$ is a rational function. At the infinitely many zeros of $f''$ (and therefore $b$-points of $f'$ and $f$) $h$ has the constant value $nb/(a-b)$. It follows 
$$
\frac{(n+1)f'' - nf'}{f-a} \equiv \frac{nb}{a-b}.
$$
This linear differential equation can be solved explicitly. With
\begin{align} \label{lambda1_lambda2}
\lambda_1 := \frac{n}{2(n+1)}, \quad \lambda_2 := \sqrt{ \frac{n^2}{4(n+1)^2}+\frac{n}{n+1} \cdot \frac{b}{a-b}}
\end{align} 
it follows that $f$ has the form
$$
f(z) = c_1 \exp[(\lambda_1 + \lambda_2)z] + c_2 \exp[(\lambda_1 - \lambda_2)z] + a.
$$
There is no loss of generality if we assume that $z=0$ is an $a$-point of $f$. This implies $c_1=-c_2$. Hence $f$ can be expressed
$$
f(z) = c \exp(\lambda_1 z) \sinh(\lambda_2 z) + a.
$$
The $a$-points $z_j$ of $f$ are exactly the zeros of $\sinh(\lambda_2 z)$ and therefore $z_j = i j \pi / \lambda_2$ with $j \in \Z$. Taking the derivative gives
$$
f'(z) = c \lambda_1 \exp(\lambda_1 z) \sinh(\lambda_2 z) + c \lambda_2 \exp(\lambda_1 z) \cosh(\lambda_2 z).
$$
With $z=0$ it follows $\lambda_2 = a/c$ and thus
$$
f'(z_j) = a \exp(\lambda_1 z_j) (-1)^j.
$$
The equation $f'(z_j)=a$ is therefore equivalent to 
$$
\exp \left( \frac{\lambda_1}{\lambda_2} i j \pi \right) = (-1)^j.
$$
It follows that $\lambda_1 = m \lambda_2$ with an odd integer $m$. But this is not compatible with (\ref{lambda1_lambda2}) and (\ref{relation_a_b}): Taking the square of
$$
\frac{n}{2(n+1)} = m \sqrt{ \frac{n^2}{4(n+1)^2}+\frac{n}{n+1} \cdot \frac{b}{a-b}}
$$
and using (\ref{relation_a_b}) gives after simple calculations
\begin{align} \label{m_n_k}
\left(1 - \frac{1}{m^2} \right)\left[(n+1)^2(k+1)+n\right]  = 4 (n+1).
\end{align}
For $m=1$ this is clearly impossible. Otherwise, since $m$ is odd, we have $1 - 1/m^2 \ge 8/9$ and already
$$
\frac{8}{9} (n+1)^2(k+1) > 4 (n+1)
$$ 
for $n \ge 2$ and $k \ge 1$ shows that (\ref{m_n_k}) can not hold. We arrive at a contradiction. Hence there are at least two $a$-points where $f''$ has the same value and $C_{n+1} \neq 0$ for all $n \ge 2$ and thus $f$ is periodic. 
$\square$

\section{Polynomials in $\exp(\lambda z)$}

\begin{rem} \label{fundamental_domain} \rm For a periodic solution of (P) let $T$ be a period with minimal absolute value, a so called {\it fundamental period}. Since a transcendental entire function is never doubly periodic, we have exactly two fundamental periods, namely $\pm T$. Since $f$ is periodic, its behaviour is completely determined on a {\it fundamental domain}, which can be chosen to be an infinite strip of width $|T|$ and such that the boundary lines are orthogonal to $\pm T$. By identifying the boundary lines in the usual way we obtain a Riemann surface $F$ which we continue to call {\it the fundamental domain} of $f$.
\end{rem}

\begin{lemma} \label{P(exp)} Let $f$ be a transcendental entire function that solves (P) with fundamental periods $\pm T$. If we set $\lambda := 2 \pi i / \pm T$, where the sign in $\pm T$ has to be chosen in an appropriate manner, then $f$ is of the form 
$$
f(z) = P \left( \text{\rm e}^{\lambda z} \right)
$$
with a unique polynomial $P$. 
\end{lemma}

\noindent
\pr By Theorem \ref{f_periodic} $f$ is periodic, and it is a well known fact (see e.g.~\cite{A}, chapter 7, section 1.1) that $f$ must be of the form $f(z) = P(\exp(\lambda z))$ with a function $P$ that is holomorphic in $\C \setminus \{ 0\}$. Since $\exp(\lambda z)$ is locally biholomorphic it follows that $P$ is uniquely determined up to the choice of the sign in $\pm T$.

We show that $P$ has to be rational. Assume to the contrary that $\infty$ is an essential singularity of $P$. (If $0$ is an essential singularity of $P$ consider $-\lambda$ instead of $\lambda$.) Since $P$ is holomorphic in a neighbourhood of $\infty$ we conclude by the characterization of Julia exceptional functions given in \cite{Os} (see also \cite{LeVi}) that there exists a sequence $w_n \to \infty$ with $|w_n| P^\#(w_n) \to \infty$ as $n \to \infty$. Choose $z_n \to \infty$ such that $w_n = \exp(\lambda z_n)$. Then $f^\#(z_n) = |\lambda \exp(\lambda z_n)| P^\#(\exp(\lambda z_n)) = |\lambda w_n| P^\#(w_n) \to \infty$ in contradiction to the boundedness of $f^\#$, which was proven in Theorem \ref{exponential_type}.

Now that we know that that $P(t)$ is rational,  proving that it is a polynomial is elementary.
If not, we can write 
$$f(z)=a+\frac{Q(t)}{t^m}$$
with $t=\exp(\lambda z)$ and a polynomial $Q$ with $Q(0)\neq 0$.  Without loss of generality 
(replacing $t$ by $1/t$ if necessary) we can assume $0<m\leq \deg(Q)$. Then
$$f'(z)=\frac{-\lambda mQ(t)}{t^m}+\frac{\lambda Q'(t)}{t^{m-1}}.$$
So if $f(\text{\rm e}^{\lambda z_0})=a$, that is $Q(t_0)=0$, then $f=a\Rightarrow f'=a$ implies
$\frac{\lambda Q'(t_0)}{t_0^{m-1}}=a$. This shows that $Q(t)$ has $\deg(Q)$ different simple zeros, and 
leads to a contradiction because the rational function $\frac{\lambda Q'(t)}{t^{m-1}}$ cannot take the value 
$a$ more than $\deg(Q)-1$ times.
$\square$

\begin{defi} \label{d_j_k} \rm 
Let $f$ be a transcendental solution of (P) such that $f''$ has zeros. Let $P$ be the above unique polynomial 
and $\lambda \neq 0$ such that $f(z) = P(\exp(\lambda z))$. We denote the degree of $P$ by $d$
and as before by $k$ the order of the zeros of $f''$, which is the same for every zero of $f''$ by Lemma \ref{g=(k+1)b/(a-b)}.

We set $t := \exp(\lambda z)$ and therefore $f(z) = P(t)$. Since $f'(z) = \lambda t P'(t)$ we define 
$\partial_z$ to be the derivation $\partial_z Q(t) := \lambda t Q'(t)$ for every function $Q$ that is 
holomorphic in $\C \setminus \{ 0 \}$. 

Instead of the usual counting functions from Nevanlinna Theory one can now use much simpler counting 
functions. We write $f_F$ whenever we consider $f$ as a function restricted to its fundamental domain $F$.
For $c\in\C$ denote by
$$n(c,f_F)\ \ \hbox{\rm and}\ \ \overline{n}(c,f_F)$$
the number of $c$-points of $f$ in the fundamental domain $F$, counted with resp. without their multiplicities.
Since $\exp(\lambda z)$ maps $F$ biholomorphically onto $\C \setminus \{ 0 \}$, these numbers simply equal the 
number of $c$-points of $P$ on $\C \setminus \{ 0 \}$, counted with resp. without multiplicities. 
Likewise $n(c,f'_F)$ equals the number of $c$-points of $\partial_z P$ on $\C \setminus \{ 0 \}$ and so on.
\end{defi}

\vskip.3cm
The basic relations between these quantities are as follows.

\begin{lemma}\label{countingfct}
Let $f$ satisfy the conditions of Definition \ref{d_j_k}. Then
\begin{itemize}
\item[(a)]  All $\partial_z^n P(t)$ have degree $d$ as polynomials in $t$.
\item[(b)]  $P(t)$ has an $a$-point at $t=0$ with multiplicity $j$ where $1\leq j <d$.
\item[(c)]  $\overline{n}(a,f_F)=n(a,f_F)=d-j$.
\item[(d)]  All $\partial_z^n P(t)$ with $n\geq 1$ have a zero of multiplicity $j$ at $t=0$.
\item[(e)]  $n(0,f''_F)=d-j$.
\item[(f)]  $n(b,f'_F)=d$ and $\overline{n}(b,f'_F)=j$.
\end{itemize}

\end{lemma}

\noindent
\pr
(a) is easy by induction.
(b) If $P(0)\neq a$, then $f$ has $d$ different $a$-points in $F$, which by the conditions of the uniqueness problem all are zeros of $f' -f$. But since the degree of $\partial_z P(t)-P(t)$ is at most $d$, then there cannot
be any zeros of $f' -f$ that are $b$-points of $f'$, in contradiction to Lemma \ref{derivative_finitely_many_B}.
On the other hand, $j=d$ is also not possible, because then $f$ would have no $a$-points at all. According to Lemma \ref{f_finitely_many_A} it would follow that $f''$ had no zeros, contradicting the assumptions.
(c) follows from (b). 
(d) is straightforward from $P(t)=a_d t^d+\cdots +a_j t^j +a$ by induction.
(e) follows from (a) and (d).
The first equation in (f) is clear because $\partial_z P(0)=0\neq b$. The second equation follows by combining 
the first with (e), keeping in mind that by Lemma \ref{g=(k+1)b/(a-b)} every zero of $f''$ is a $b$-point of $f'$.
$\square$

\vskip.3cm
The following estimates are less obvious.

\begin{lemma} \label{d_and_j} Let $f$, $P$, $\lambda$ and $d$ fulfil the assumptions of Definition \ref{d_j_k} and let $j$ be as in Lemma \ref{countingfct}. Then the following inequality holds:
\begin{equation} \label{d_und_j}
2j \le d \le 3j.
\end{equation}
\end{lemma}

\noindent
\pr First we consider
\begin{equation} \label{h1}
h_1 := \frac{f'(z)}{f(z)-a} = \frac{\partial_z \, P(t)}{P(t)-a} = \frac{\lambda t P'(t)}{P(t)-a}.
\end{equation}
As a function of $t$ the degree of $h_1$ is $d-j$. This follows from the fact that neither $t=0$ nor $t=\infty$ is a pole of $h_1$: In $t=0$ the function $P(t)-a$ has a zero of order $j$, and so does $t P'(t)$ and since $d$ is the degree of both $P(t)-a$ and $t P'(t)$ there is no pole in $t=\infty$. So the poles of $h_1$ are the $a$-points of $P$ in $\C \setminus \{ 0 \}$, whose number is $d-j$ by Lemma \ref{countingfct} (c).
On the other hand again from Lemma \ref{countingfct} we know  
that $j=\overline{n}(b,f'_F)$ is a lower bound for the number of $b/(b-a)$-points of $h_1$. We get $j \le d-j$ which is the first inequality we have to prove.

For the second inequality we consider
\begin{equation} \label{h2}
h_2 := \frac{f''(z)}{f'(z)-b} = \frac{\partial_z^2 P(t)}{\partial_z P(t)-b} = \frac{\lambda^2 t P'(t) + \lambda^2 t^2 P''(t)}{\lambda t P'(t)-b}.
\end{equation}
and argue in the same manner. $h_2$ has no pole in $t=0$ since $\lambda t P'(t)-b$ has the value $-b \neq 0$ at $t=0$. Further $h_2$ has no pole in $t=\infty$ because the degree of $t P'(t) + t^2 P''(t)$ is less than or equal to the degree of $\lambda t P'(t)-b$.  The poles of $h_2$ are exactly the $b$-points of $f'$ and as poles of $h_2$ they are all simple because $h_2$ is a logarithmic derivative. So the degree of $h_2$ is $\overline{n}(b,f'_F)=j$.
We consider the $d-j$ different $a$-points of $P$ in $\C \setminus \{ 0 \}$ (see Lemma \ref{countingfct} (c)). At these we have by assumption $\lambda t P'(t) = a$ and $f'' = \lambda^2 t P'(t) + \lambda^2 t^2 P''(t)$ takes by Lemma \ref{f''_at_a,a} one of two possible values, say $u$ or $v$. It follows that the $a$-points of $P$ in $\C \setminus \{ 0 \}$ are contained in the $u/(a-b)$-points and $v/(a-b)$-points of $h_2$ and therefore $d-j \le 2j$ which completes the proof. 
$\square$

\section{The $b$-points of $f'$}\label{sectionbpoints}

\begin{lemma} \label{f''_at_simple_b} Let $f$ solve (P). If $f'$ has a simple $b$-point at $z_0$ and if $f''$ has zeros of order $k \ge 1$ (elsewhere), then $f(z_0) = b$ and $f''(z_0) = -kb$
\end{lemma}

\noindent
\pr From Lemma \ref{g=(k+1)b/(a-b)} we know $g \equiv (k+1)b/(a-b)$ and $f(z_0) = f'(z_0) = b$. This implies
$$
\frac{f'-f}{f'-b} = \frac{(k+1)b}{a-b} \frac{f-a}{f''}.
$$
Taking the limit for $z \to z_0$ on both sides, we obtain
$$
\frac{f''(z_0)- b}{f''(z_0)} = \frac{-(k+1)b}{f''(z_0)},
$$
which proves the claim. 
$\square$

\begin{lemma} \label{unique_Taylor_in_F} Let $f$ be a transcendental entire function that solves (P) and $z_1$ and $z_2$ be two points in the fundamental domain $F$ of $f$ with $f(z_1) = f(z_2)$. Then the Taylor coefficients of $f$ at $z_1$ and $z_2$ do not coincide.
\end{lemma}

\noindent
\pr If the Taylor coefficients coincide, then from Lemma \ref{f_periodic_taylor} we know that $t = z_1 - z_2$ is a period of $f$. If $t$ and $\pm T$ were collinear, then $\pm T$ would not be minimal. Otherwise $t$ is a period that is independent of $T$, therefore $f$ would be doubly periodic, a contradiction. 
$\square $

\begin{prop} \label{simple_b} If $f$ is a transcendental entire function that solves (P) such that $f''$ has zeros, then $f'$ has at most one simple $b$-point in the fundamental domain.
\end{prop}

\noindent
\pr By Remark \ref{fundamental_domain} we have to prove that at every simple $b$-point of $f'$ the Taylor coefficients of $f$ are the same. Suppose that $z$ is a simple $b$-point of $f'$. By Lemma \ref{f''_at_simple_b} we have $f''(z) = -kb$ where $k$ is again the order of the zeros of $f''$. Exactly as in the proof of Theorem \ref{f_periodic} we get the equation
\begin{align*}
f^{(n+2)}(f'-f) + {n \choose 1} f^{(n+1)}(f''-f') + \ldots + f''(f^{(n+1)} - f^{(n)}) & \\
= \frac{(k+1)b}{a-b}\left[f^{(n)}(f'-b) + \ldots + (f-a) f^{(n+1)} \right] &.
\end{align*}
By assumption we have $f'(z) - f(z) = 0$, and therefore the first term on the left disappears. Again we collect the terms in which $f^{(n+1)}$ occurs on the left side of the equation. Since the last term on the right does not vanish in general, the coefficient of $f^{(n+1)}$ in this case is
$$
\widetilde{C}_{n+1}:= (n+1)f'' - nf' -   \frac{(k+1)b}{a-b} (f-a).
$$
Using $f=f'=b$ and $f''=-kb$ we get $\widetilde{C}_{n+1} = [1-n(k+1)]b$ which is never $0$ for $n \ge 1$ and $k \ge 1$. Therefore the value of $f^{(n+1)}$ at $z$ can be determined uniquely for every $n \ge 2$ (not depending on $z$). 
$\square$
\vskip.3cm

Inspired by this result we will prove our main theorem with a case distinction depending on whether $f'$ has simple $b$-points or not. We use the notation from Definition \ref{d_j_k} and Lemma \ref{countingfct} for the rest of the section.

\begin{lemma} \label{b_points_multiple} If $f$ solves (P) and all $b$-points of $f'$ are multiple then $j=1$. Furthermore, then either $d=3$ (and $k=2$) or $d=2$ (and $k=1$).
\end{lemma}

\noindent
\pr From the previous results, in particular Lemmas \ref{P(exp)},  \ref{countingfct} and \ref{g=(k+1)b/(a-b)} (c) 
we have
$$f'(z)=C\prod_{s=1}^{j}(t-r_s)^{k+1}+b$$
with $t=\text{\rm e}^{\lambda z}$ where 
$r_1,\ldots, r_s$ are nonzero and pairwise distinct.
\par
Writing $Q=\prod_{s=1}^{j}(t-r_s)$ we have
$$f' -b =CQ^{k+1}$$
and also
$$f'=c_d t^d +\cdots + c_j t^j.$$
Let $w$ be the smallest positive integer for which $Q$ has a nonzero term $q_w t^w$. Then
$$Q^{k+1}=
(q_0 +q_w t^w +\cdots +q_j t^j)^{k+1}=
q_0^{k+1}+(k+1)q_0^k q_w t^w +\cdots$$
So $w\geq j$, and since $Q$ actually has degree $j$, it must be of the form $Q=q_j t^j +q_0$.
Thus $f$ is actually a polynomial in $t^j$. So $\frac{2\pi i}{j\lambda}$ is a period.  But since 
$\frac{2\pi i}{\lambda}$ 
is a fundamental period,  we must have $j=1$. With this Lemma \ref{d_and_j} implies $d\leq 3$.
$\square$

\begin{lemma} \label{d=3_j=1_k=2} There are no solutions of (P) with $d=3$, $j=1$ and $k=2$.
\end{lemma}

\noindent
\pr From 
$$f=Ct(t-r)(t-s)+a$$
with $t=\text{\rm e}^{\lambda z}$ we obtain
$$f'=\lambda Ct(3t^2 -2(r+s)t+rs).$$
Since $t=r$ and $t=s$ are $a$-points of $f$, and hence must also be $a$-points of $f'$ we have
$$a=\lambda Cr^2 (r-s)$$
and 
$$a=\lambda Cs^2 (s-r).$$
Subtraction gives 
$0=\lambda C(r-s)(r^2 +s^2)$, and hence $s=\pm ir$. (Note that $r=s$ is not possible because the $a$-points of $f$ are simple.) With this one easily checks that the discriminant of the quadratic factor in
$$f''=\lambda^2 Ct(9t^2 -4(r+s)t+rs)$$
is nonzero, so $f''$ does not have a double root. 
$\square$

\begin{prop} \label{iv)} Solutions $f$ of (P) for which all $b$-points of $f$ are multiple exist if and only if $b=\frac{-a}{8}$. They are exactly the solutions of type (iv) in Theorem \ref{main}, that is,
$$ f(z) = 
6aC\text{\rm e}^{\frac{z}{6}}(C\text{\rm e}^{\frac{z}{6}}-1)+a$$
with a nonzero constant $C$.
\end{prop}

\noindent
\pr From Lemma \ref{b_points_multiple} and Lemma \ref{d=3_j=1_k=2} we know that only $d=2$ and $k=1$ is possible. Thus $f$ must have the form
$$f=At(t-r)+a \ \ \hbox{\rm with}\ \ t=\text{\rm e}^{\lambda z},$$
and hence
$$f' = \lambda At(2t-r)\ \ \ \hbox{\rm and}\ \ f''=\lambda^2 At(4t-r).$$
Since $t=r$ (the $a$-point of $f$) must also be an $a$-point of $f'$, we get $a=\lambda Ar^2$. 
\par
Moreover, the zero of $f''$, that is $t=\frac{r}{4}$, must also give the double $b$-point of $f'$, whence 
$b=\frac{-\lambda Ar^2}{8}$, and hence $b=\frac{-a}{8}$.
\par
Since $t=\frac{r}{4}$ must also be a $b$-point of $f$, we have 
$b=A\frac{r}{4}\frac{-3r}{4}+a$, which yields $\lambda=\frac{1}{6}$.
\par
Setting $C=\frac{1}{r}$ gives the desired form. 
$\square$

\section{The end of the proof}\label{sectionendofproof}

This last section has the purpose to show that (P) has no solution where $f'$ has simple and multiple $b$-points. It follows that either $f''$ has no zeros, which gives the solutions (ii) and (iii) in Theorem \ref{main}, or all $b$-points of $f'$ are of order $2$ which according to Proposition \ref{iv)} leads to the functions in (iv) of Theorem \ref{main}.

\begin{lemma} \label{relations} Let $f$, $P$, $\lambda$, $d$ and $k$ fulfil the assumptions of Definition \ref{d_j_k} and let $j$ be as in Lemma \ref{countingfct}. Then the following relation holds:
\begin{align}
\frac{(k+1)b}{a-b} & = \frac{d^2 j^2 a b}{(d^2 b - j^2 a)^2}. \label{ohne_lambda}
\end{align}   
\end{lemma}

\noindent
\pr We can apply Lemma \ref{g=(k+1)b/(a-b)} and therefore
$$
g(z)= \frac{f''(z)(f'(z)-f(z))}{(f(z)-a)(f'(z)-b)} \equiv \frac{(k+1)b}{a-b}.
$$
This translates for $P(t)$ to
$$ \label{g_for_P}
g(t)= \frac{(\lambda^2 t P'(t) + \lambda^2 t^2 P''(t))(\lambda t P'(t)- P(t))}{(P(t)-a)(\lambda t P'(t)-b)} \equiv \frac{(k+1)b}{a-b},
$$
so
$$ \label{g_for_P_2}
(\lambda^2 t P'(t) + \lambda^2 t^2 P''(t))(\lambda t P'(t)- P(t)) = \frac{(k+1)b}{a-b} (P(t)-a)(\lambda t P'(t)-b).
$$
Setting
$$
P(t) = a + a_j t^j + \ldots + a_d t^d.
$$
and comparing the highest and lowest terms on both sides yields the equations 
\begin{align}
\frac{(k+1)b}{a-b} & = (\lambda d)^2 - \lambda d \label{lambda_d}
\end{align}
and \begin{align}
\frac{(k+1)b}{a-b} & = \frac{a}{b} (\lambda j)^2. \label{lambda_j}
\end{align} 
It follows that the right hand sides of (\ref{lambda_d}) and (\ref{lambda_j}) are equal. After cancelling $\lambda \neq 0$ this gives
$$
\lambda d^2 - d = \frac{a}{b} \lambda j^2,
$$
thus
$$
\lambda = \frac{d b}{d^2 b - j^2 a}.
$$
Substituting this back into (\ref{lambda_j}) gives (\ref{ohne_lambda}). 
$\square$

\begin{lemma} \label{k_n_square} Let $f$, $P$, $\lambda$, $d$ and $k$ fulfil the assumptions of Definition \ref{d_j_k}, let $j$ be as in Lemma \ref{countingfct} and suppose that $C_{n+1}=0$ for some $n \ge 2$ at an $a$-point of $f$, with $C_{n+1}$ from (\ref{C_{n+1}}). Then the following equality holds
\begin{align} \label{d_j_equation}
d^2 j^2 (n+1)^2 \left( (k+1)(n+1)^2 + n \right) = \left(d^2 n + j^2 (n+1)^2 (k+1)\right)^2.
\end{align}
In particular $(k+1)(n+1)^2 + n$ has to be a square.
\end{lemma}

\noindent
\pr Since $C_{n+1}=0$ at an $a$-point, we have according to (\ref{relation_a_b}) the relation $a = -(n+1)^2 (k+1)b/n$. Using this in (\ref{ohne_lambda}) gives after some simple cancelations
\begin{align}
\frac{1}{(k+1)(n+1)^2/n+1} & = \frac{d^2 j^2(n+1)^2/n }{(d^2  + j^2(n+1)^2 (k+1)/n)^2}.
\end{align}
Now (\ref{d_j_equation}) follows easily. Since $d^2 j^2 (n+1)^2$ and $\left(d^2 n + j^2 (n+1)^2 (k+1)  \right)^2$ are squares, $(k+1)(n+1)^2 + n$ also has to be a square. 
$\square$

\begin{lemma} \label{2<=k<=4}
If $f$ solves (P) and $f'$ has simple and multiple $b$-points, then 
$$2\leq k\leq 4.$$
\end{lemma}

\noindent
\pr  
From Lemma \ref{countingfct} (f) we know $\overline{n}(b,f'_F)= j$. By Proposition \ref{simple_b} there can be at most one simple $b$-point for $f'_F$. Hence the number of ramified $b$-points is $j-1$, which is not zero by assumption. The common multiplicity then is $k = (d-j)/(j-1) \leq 2j/(j-1)=2 + 2/(j-1) \le 4$ by Lemma \ref{d_and_j}. On the other hand, again by Lemma \ref{d_and_j}, we have 
$k=\frac{d-j}{j-1}\geq \frac{2j-j}{j-1}>1$. 
$\square$
\vskip.3cm
We now show that for the only possible $k$ with $2\leq k\leq 4$ the expression $(k+1)(n+1)^2 + n$ is never a square for positive $n$.

\begin{lemma} \label{k=3,4}
\begin{itemize}
\item[(a)]  $5(n+1)^2 +n$ is not a square for any positive integer $n$.
\item[(b)]  $4(n+1)^2 +n$ is not a square for any positive integer $n$.
\end{itemize}
\end{lemma}

\noindent
\pr (a)   Plugging in 
$n=0,\ 1,\ 2,\ 3,\ 4,\ 5,\ 6,\ 7,\ 8$ we see that for any integer $n$ the expression 
$5(n+1)^2 +n$ is congruent to $2$, $3$, $5$ or $8$ modulo $9$. But the squares modulo $9$ are congruent to $0$, $1$, $4$ or $7$ modulo $9$.
\par
(b)  We want to solve $4(n+1)^2 +n=s^2$ in positive integers $n$ and $s$. Multiplying this equation with $16$ and setting $x=8n+9$, $y=4s$ we get $x^2 -y^2 =17$. Since $17=(x+y)(x-y)$ must be a factorization of $17$ (with positive integers) only $x+y=17$, $x-y=1$ is possible. So $x=9$, $y=8$ and $n=0$ is the only solution. 
$\square$

\begin{lemma} \label{k=3}
$3(n+1)^2 +n$ is not a square for any positive integer $n$.
\end{lemma}

\noindent
\pr Since for $n=-11$ this gives a square, it is not possible to prove this lemma by congruence calculation.
\par
Multiplying the equation 
$3(n+1)^2 +n=s^2$ with $12$ and setting $x=6n+7$ and $y=2s$ we obtain
$$x^2 -3y^2 =13.$$
There are solutions of this equation, for example $(x,y)=(4,1)$ or $(5,2)$ or $(-5,2)$. But we need solutions where $y$ is even and $x$ is positive and congruent to $1$ modulo $6$.
\par
Motivated by 
$(x+y\sqrt{3})(x-y\sqrt{3})=13$ as in the previous proof, the key idea is to associate with each solution $(x,y)$ the real number $x+y\sqrt{3}$. Let $x$ and $y$ be positive. Then obviously $x+y\sqrt{3}>\sqrt{13}$ and $\sqrt{13}>x-y\sqrt{3}>0$, whereas $-x\pm y\sqrt{3}$ both are negative.
\par
The key step is to show that if there is a solution $(x,y)$ with $y$ positive and even and $x$ positive and congruent to $1$ modulo $6$, then there also is a solution $(x_n,y_n)$ with these properties and
$$\sqrt{13}<x_n +y_n \sqrt{3}\leq 51.$$
If $x+y\sqrt{3}>51$ we divide it by $7+4\sqrt{3}\approx 13.92$ or in other words multiply it by $7-4\sqrt{3}$. The result is $x_1 +y_1\sqrt{3}$ with 
$x_1 =7x-12y$ and $y_1 =7y-4x$. 
One easily checks that $(x_1,y_1)$ is again a solution with $y_1$ being even and $x_1$ congruent to $1$ modulo $6$. But 
$x_1 +y_1\sqrt{3}$ is almost $14$ times smaller than $x+y\sqrt{3}$. 
\par
Continuing like this we eventually get a solution $(x_n,y_n)$ with
$51\geq x_n +y_n\sqrt{3}>\sqrt{13}$ (because $51/14>\sqrt{13}$). By the discussion above this also shows that $x_n$ and $y_n$ are positive.
\par
Now we can simply plug in the finitely many pairs of positive integers $(x,y)$ with $x+y\sqrt{3}\leq 51$, $x$ congruent to $1$ modulo $6$ and even $y$. None of these solve $x^2 -3y^2 =13$, hence there is no solution $(x,y)$. 
$\square$

\begin{rem} \label{Remark_after_square} \rm a) Although the proof of Lemma \ref{k=3} is selfcontained and understandable on an elementary level, we do not want to hide the structural ideas behind it. All of them are applications of fundamental theorems from algebraic number theory to the ring
$$R=\{x+y\sqrt{3}\ :\ x,y\in\Z\}.$$
The first one is that the map 
$N\ :\ x+y\sqrt{3}\mapsto 
(x+y\sqrt{3})(x-y\sqrt{3})=x^2 - 3y^2$ 
from $R$ to $\Z$ (called the norm map) is multiplicative. From this it follows easily that the units of $R$ are exactly the elements with norm $\pm 1$, and that if $N(r)=\pm p$ for a prime $p$ then $r$ is irreducible. Moreover, $R$ is a Euclidean ring with respect to $N$, and hence has unique prime factorization. From this and $13=(4+\sqrt{3})(4-\sqrt{3})$ it follows that every element $x+y\sqrt{3}$ with $(x+y\sqrt{3})(x-y\sqrt{3})=13$ differs multiplicatively from either $4+\sqrt{3}$ or from $4-\sqrt{3}$ by a unit of $R$.
\par
By another fundamental theorem of algebraic number theory $\varepsilon=2+\sqrt{3}$ is the fundamental unit of $R$. This means that every unit of $R$ is of the form 
$\pm \varepsilon^n$ with $n\in\Z$. Since $\varepsilon=2+\sqrt{3}\approx 3.73$, by multiplying or dividing by units we can get solutions $(x,y)$ with $x+y\sqrt{3}$ as big or as small as we want. But multiplying (or dividing) $x+y\sqrt{3}$ by $\varepsilon$ switches the parities of $x$ and $y$. This is the reason why in the proof we divide by
$\varepsilon^2 =7+4\sqrt{3}$.

b) Note that the Lemmata \ref{k_n_square}, \ref{2<=k<=4}, \ref{k=3,4} and \ref{k=3} show that $C_{n+1} \neq 0$ for all $n \ge 2$ if $f'$ has simple and multiple $b$-points.
\end{rem}

For the rest of the section we again use the notation from Definition \ref{d_j_k} and Lemma \ref{countingfct}.

\begin{lemma} \label{d=4_j=2_k=2}
If $f$ is a solution of (P) such that $f'$ has multiple and simple $b$-points, then
$d=4$, $j=2$ and $k=2$.
\end{lemma}

\noindent
\pr By Lemma \ref{f''_at_a,a} at every $a$-point of $f$ there are only two different possible values for $f''$.  By Remark \ref{Remark_after_square}, b) we know that $C_{n+1} \neq 0$ for all $n \ge 2$, thus $f^{(n+1)}$ is uniquely determined at every $a$-point, depending only on the value of $f''$. Hence there are only two possible sequences of Taylor coefficients at $a$-points, and by Lemma \ref{unique_Taylor_in_F} it follows that $f_F$ has at most two $a$-points, hence $d-j \le 2$. By Lemma \ref{d_and_j} it follows $j \le d - j \le 2$ and therefore $d = d - j + j \le 4$. On the other hand, since $f'$ has at least two $b$-points (a simple one and a ramified one) we get from Lemma \ref{countingfct} (f) that $j\geq 2$ and thus $j=2$. Lemma \ref{d_and_j} gives $d\geq 4$, so $d=4$. It follows $k = (d - j)/(j-1) = 2$. 
$\square$

\begin{lemma}
There are no solutions of (P) with $d=4$, $j=2$ and $k=2$, and hence no solutions of (P) such that $f'$ has multiple and simple $b$-points.
\end{lemma}

\noindent
\pr This is practically the same proof as for Lemma \ref{d=3_j=1_k=2}. With the Ansatz
$$f=Ct^2 (t-r)(t-s)+a$$
we obtain from $f=a\Rightarrow f'=a$ and
$$
f'= \lambda Ct^2(4t^2 -3(r+s)t+2rs)
$$
that $s/r$ is a $6$-th root of unity different from $1$. And with that we see that
$$f''= \lambda^2 Ct^2(16t^2 -9(r+s)t+4rs)$$
does not have the required double root. The last assertion follows directly from Lemma \ref{d=4_j=2_k=2}. 
$\square$


\vskip.5cm
\noindent
Andreas Sauer, \\
Hochschule Ruhr West, \\
Duisburger Str.~100, \\
45479 M\"{u}lheim an der Ruhr,\\
Germany\\
andreas.sauer@hs-ruhrwest.de
\\ \\
Andreas Schweizer, \\
Am Felsenkeller 61,\\
78713 Schramberg,\\
Germany\\

\end{document}